\documentclass[11pt,a4paper]{article}
\usepackage[T1]{fontenc}
\usepackage[utf8]{inputenc}
\usepackage{lmodern,microtype}
\usepackage[margin=1in]{geometry}
\usepackage{amsmath,amssymb,amsthm,mathtools}
\usepackage{booktabs,enumitem}
\usepackage[hidelinks]{hyperref}
\hypersetup{pdftitle={Second-Order Threshold and Value Asymptotics in Multiple-Stopping Secretary Problems},pdfauthor={José María Grau Ribas; Raúl Sánchez Galán}}
\theoremstyle{plain}
\newtheorem{theorem}{Theorem}[section]
\newtheorem{proposition}[theorem]{Proposition}
\newtheorem{lemma}[theorem]{Lemma}
\newtheorem{corollary}[theorem]{Corollary}
\theoremstyle{remark}
\newtheorem{remark}[theorem]{Remark}
\newcommand{\e}{\mathrm e}
\newcommand{\N}{\mathbb N}

\newcommand{\R}{\mathbb R}
\newcommand{\E}{\mathbb E}
\newcommand{\Prob}{\mathbb P}
\newcommand{\floor}[1]{\left\lfloor #1\right\rfloor}
\newcommand{\fract}[1]{\left\{#1\right\}}
\newcommand{\esym}{\mathsf e}
\newcommand{\lam}{\boldsymbol\lambda}
\newcommand{\doi}[1]{\href{https://doi.org/#1}{\nolinkurl{doi:#1}}}
\numberwithin{equation}{section}

\title{Second-Order Threshold and Value Asymptotics\\
in Multiple-Stopping Secretary Problems}
\author{Jos\'e Mar\'ia Grau Ribas \and Ra\'ul S\'anchez Gal\'an}
\date{}

\begin{document}
\maketitle

\begin{abstract}
We study the classical secretary problem with a fixed number of selection
opportunities. We derive explicit corrections to the limiting thresholds
and to the optimal success probability. A cancellation identity removes
the influence of later threshold rounding from the first correction,
while a Poisson representation explains its probabilistic meaning.
Every bounded rounding of the limiting thresholds shares the first value
correction. We obtain integer rounding laws, eventual threshold identities
along continued-fraction convergents, and joint rounding frequencies.
Rigorous rational bounds show that the first eleven corrections are
positive and the twelfth is negative, producing the first shifted
convergent identity.
\end{abstract}



\section{Introduction}\label{sec:intro}

In the classical secretary problem, $N$ candidates arrive in a uniformly
random order. Only their relative ranks are observed, decisions are
irrevocable, and success means selecting the best candidate. The familiar
single-choice rule rejects an initial proportion $1/\e$ and then selects
the next record (a candidate better than all predecessors). With $m$
selection opportunities, success means that at
least one of the selected candidates is the overall best. Optimal rules
again have a threshold structure, now with one threshold for each remaining
opportunity. The classical development is described by Gilbert and
Mosteller~\cite{GilbertMosteller1966} and
Ferguson~\cite{Ferguson1989}.

Write $i_N^{(k)}$ for the first eligible index when $k$ opportunities
remain, and use rejected-count coordinates
\[
 r_N^{(k)}=i_N^{(k)}-1.
\]
The thresholds can be chosen consistently as $m$ varies, with
$i_N^{(m)}\le\cdots\le i_N^{(1)}$. The record-indicator formulation turns
the problem into selection of the last success in independent Bernoulli
trials. Bruss' odds theorem~\cite{Bruss2000} covers one selection, and the
multiple-selection theorem of Ano, Kakinuma, and
Miyoshi~\cite{AnoKakinumaMiyoshi2010} establishes the nested optimal rule.
Matsui and Ano~\cite{MatsuiAno2016} identify positive rational constants
$\lambda_j$ such that, for each fixed $k$,
\begin{equation}\label{eq:intro-leading}
 \frac{i_N^{(k)}}N\longrightarrow\ell_k,
 \qquad
 \ell_k=\exp\left(-\sum_{j=1}^k\lambda_j\right).
\end{equation}
They also identify the limiting success probability as
$\sum_{k=1}^m\ell_k$. A complementary continuous linear-programming
approach of Chan, Chen, and Jiang~\cite{ChanChenJiang2015} gives the
rational-exponent description and the same value identity for the
$(m,1)$-secretary problem.

The proportions in \eqref{eq:intro-leading} do not specify which integer
threshold is optimal. There is a further difficulty in the multiple-choice
problem: the payoff near the $k$th threshold depends on the later integer
thresholds, whose rounding errors need not converge. Our first result
shows that these errors cancel at the first correction order. If
$\Delta_{k,N}(r)$ is the increment in success probability when the $k$th
rejected-count threshold is changed from $r-1$ to $r$, with the later
thresholds fixed optimally, then
\begin{equation}\label{eq:intro-local}
 N^2\Delta_{k,N}(\ell_kN+u)
 =-\frac{u-c_k}{\ell_k}+O(N^{-1})
\end{equation}
uniformly for bounded admissible offsets $u$. The constant $c_k$ is an
explicit combination of diagonal second derivatives of the winning-pattern
polynomial; see \eqref{eq:c-definition}. The cancellation follows from
the identity that every first derivative of that polynomial equals one
at the limiting block parameters.

Our second result is the finite-horizon refinement
\begin{equation}\label{eq:intro-value}
 P_N^*(m)=\sum_{k=1}^m\ell_k
       +\frac1N\sum_{k=1}^m c_k+O(N^{-2}),
\end{equation}
where $P_N^*(m)$ is the optimal success probability. In fact,
\eqref{eq:intro-value} holds for the success probability of \emph{any}
threshold vector satisfying $r_k=\ell_kN+O(1)$, uniformly over bounded
displacements. Thus arbitrary bounded rounding already attains the
optimal value to within $O(N^{-2})$, even when it does not select the
exact optimal integers. This value statement and the local cancellation
are the main contributions of the paper.

The diagonal second-derivative term has a probabilistic meaning. Independent
Poisson variables with means equal to the Bernoulli odds recover the
record indicators when conditioned to be at most one. Within a block, the
first departure from the Poisson count formula is the probability of a
collision between two points assigned to the same candidate time. This
gives an exact embedding and identifies the Hessian correction as a
weighted factorial moment conditional on a winning pattern.

The local expansion also has arithmetic consequences. For one selection,
\cite{Sanchez2026Convergents} proves that the optimal rejected-count threshold
is $p$ whenever $p/q$ is a convergent of $1/\e$ and $q\ge3$.
The scalar principle
relating affine payoff increments to integer rounding was developed in
\cite{Sanchez2026Rounding}. Here we prove the required expansion throughout
the coupled threshold hierarchy and allow corrections of either sign.
Outside a transition band of width $O(N^{-1})$,
\[
 r_N^{(k)}=\floor{\ell_kN+c_k}.
\]
Moreover, if $p_j/q_j$ are the continued-fraction convergents of $\ell_k$,
then $r_{q_j}^{(k)}=p_j+\floor{c_k}$ for all sufficiently large $j$.
We show that every $c_k$ is transcendental, so no integer correction creates
a limiting tie. The limiting joint frequencies of the rounded threshold
choices factorize, although the finite-horizon optimization is coupled. Exact finite computations establish $0<c_k<1/2$ for
$1\le k\le11$ and $-1<c_{12}<0$. Consequently the twelfth threshold is
the first one whose convergent identity has shift $-1$.

Throughout, the number of selection opportunities is fixed as $N$ tends
to infinity. No estimate uniform in a growing number of opportunities is
claimed. Section~\ref{sec:model} records the exact finite-horizon formulas.
Section~\ref{sec:expansions} proves the block expansions and derivative
identity. Sections~\ref{sec:local} and~\ref{sec:value} establish the threshold
and value results. Rounding and arithmetic consequences are proved in
Section~\ref{sec:rounding}, and Section~\ref{sec:examples} gives examples
and a reproducible finite sign certificate.

\section{The model and exact payoff formulas}\label{sec:model}

\subsection{Thresholds and winning patterns}

Let $X_t$ indicate that candidate $t$ is a record. For a random permutation,
the variables $X_2,\ldots,X_N$ are independent, with
\begin{equation}\label{eq:odds}
 p_t=\Prob(X_t=1)=\frac1t,\qquad
 q_t=1-p_t=\frac{t-1}{t},\qquad
 o_t=\frac{p_t}{q_t}=\frac1{t-1}.
\end{equation}
The best candidate is the last record. A threshold vector
$\mathbf i=(i_m,\ldots,i_1)$ satisfies
$2\le i_m\le\cdots\le i_1\le N$. Its rule selects a record whenever the
number of previous selections is smaller than the number of thresholds
already passed. Denote its success probability by $P_N(\mathbf i)$.

We work with the standard recursively compatible optimal thresholds given
by the multiple-stopping odds theorem
\cite{AnoKakinumaMiyoshi2010,MatsuiAno2016}: after fixing the optimal
$i_N^{(1)},\ldots,i_N^{(k-1)}$, maximize
\begin{equation}\label{eq:P-definition}
 P_{k,N}(i)=P_N(i,i_N^{(k-1)},\ldots,i_N^{(1)})
\end{equation}
over $2\le i\le\min\{i_N^{(k-1)},N\}$, with $i_N^{(0)}=N+1$.
We choose the largest maximizer at a tie. All statements about these
thresholds are for sufficiently large $N$ with $k$ fixed. In particular,
\eqref{eq:intro-leading} places them strictly after the first candidate;
the exceptional small horizons at which selection starts at index one
play no role.

For an arbitrary threshold vector, set $i_0=N+1$ and
\[
 B_j=\{i_j,\ldots,i_{j-1}-1\},\qquad 1\le j\le m.
\]
For a finite set $B$ and an integer $b\ge0$, let
\[
 \esym_b(B)=\sum_{\substack{A\subseteq B\\|A|=b}}\prod_{t\in A}o_t,
 \qquad \esym_0(B)=1.
\]
The value is zero when $b>|B|$. Define $\Xi_k$ to be the set of
nonnegative integer vectors $\mathbf b=(b_k,\ldots,b_1)$ for which some
$s\in\{1,\ldots,k\}$ satisfies
\begin{equation}\label{eq:patterns}
 \sum_{j=s}^k b_j=k-s+1,\qquad
 \sum_{j=t}^k b_j<k-t+1\quad(s<t\le k),\qquad
 b_{s-1}=\cdots=b_1=0.
\end{equation}
For example,
\[
 \Xi_1=\{(1)\},\quad
 \Xi_2=\{(1,0),(0,2)\},\quad
 \Xi_3=\{(1,0,0),(0,2,0),(0,1,2),(0,0,3)\}.
\]
The winning-pattern characterization of Matsui and
Ano~\cite{MatsuiAno2016} says that a record sequence wins precisely when
there is a unique $k$ for which its counts in $B_k,\ldots,B_1$ belong to
$\Xi_k$. Independence therefore gives
\begin{equation}\label{eq:exact-value}
 P_N(\mathbf i)
 =\sum_{k=1}^m\frac{i_k-1}{N}\,\mathcal H_{k,N},\qquad
 \mathcal H_{k,N}
 =\sum_{\mathbf b\in\Xi_k}\prod_{j=1}^k\esym_{b_j}(B_j).
\end{equation}
Indeed, the probability of prescribed counts in these blocks is the
product of their elementary symmetric sums times
$\prod_{t=i_k}^Nq_t=(i_k-1)/N$. The disjoint characterization justifies
summing over $k$. We will use \eqref{eq:exact-value} for both optimal and
nonoptimal threshold vectors.

Introduce the finite polynomials
\begin{equation}\label{eq:polynomials}
 H_k(x_1,\ldots,x_k)
 =\sum_{\mathbf b\in\Xi_k}\prod_{j=1}^k\frac{x_j^{b_j}}{b_j!},
 \qquad H_k=x_k+\Phi_k(x_1,\ldots,x_{k-1}).
\end{equation}
The second identity holds because the only vector in $\Xi_k$ with
$b_k>0$ is $(1,0,\ldots,0)$. The limiting parameters are the unique
positive solution of the triangular system
\begin{equation}\label{eq:lambda}
 H_k(\lambda_1,\ldots,\lambda_k)=1,
 \qquad \lambda_1=1,\qquad
 \lambda_k=1-\Phi_k(\lambda_1,\ldots,\lambda_{k-1}).
\end{equation}
Positivity and the threshold limits \eqref{eq:intro-leading} are the
first-order results of~\cite{MatsuiAno2016}. The recurrence also shows
that every $\lambda_k$ is rational. We use
\begin{equation}\label{eq:ell}
 L_0=0,\quad L_k=\sum_{j=1}^k\lambda_j,
 \qquad \ell_k=\e^{-L_k},\quad \ell_0=1,
 \qquad \delta_j=\frac1{\ell_j}-\frac1{\ell_{j-1}}>0.
\end{equation}

\subsection{An exact difference and its sign}

Fix the later optimal thresholds in \eqref{eq:P-definition}. Their blocks
$B_1,\ldots,B_{k-1}$ are now fixed, and put
\[
 \mathcal C_{k,N}
 =\sum_{\mathbf b\in\Xi_k\setminus\{(1,0,\ldots,0)\}}
           \prod_{j<k}\esym_{b_j}(B_j).
\]
For $k=1$ this is zero. Equation~\eqref{eq:exact-value} becomes
\begin{equation}\label{eq:separated-payoff}
 P_{k,N}(i)=P_N^*(k-1)
  +\frac{i-1}{N}
      \left(\sum_{t=i}^{i_N^{(k-1)}-1}\frac1{t-1}+\mathcal C_{k,N}\right),
\end{equation}
where $P_N^*(0)=0$.

\begin{lemma}[Exact difference]\label{lem:difference}
For $3\le i\le\min\{i_N^{(k-1)},N\}$,
\begin{equation}\label{eq:difference}
 P_{k,N}(i-1)-P_{k,N}(i)=\frac1N G_{k,N}(i),\qquad
 G_{k,N}(i)=1-\sum_{t=i}^{i_N^{(k-1)}-1}\frac1{t-1}-\mathcal C_{k,N}.
\end{equation}
The function $G_{k,N}$ is strictly increasing on its integer domain,
and $P_{k,N}$ is unimodal there.
\end{lemma}
\begin{proof}
Write $S(i)$ for the sum in \eqref{eq:separated-payoff}.
Since $S(i-1)=S(i)+1/(i-2)$, subtraction gives
\[
 N\bigl(P_{k,N}(i-1)-P_{k,N}(i)\bigr)
 =(i-2)(S(i)+\mathcal C_{k,N})+1
  -(i-1)(S(i)+\mathcal C_{k,N}),
\]
which is \eqref{eq:difference}. Furthermore,
$G_{k,N}(i+1)-G_{k,N}(i)=1/(i-1)>0$ whenever both sides are defined.
Thus the forward payoff increments are strictly decreasing.
\end{proof}

\section{Block expansions and the cancellation identity}\label{sec:expansions}

\subsection{Uniform expansions within a block}

\begin{lemma}[Block-odds expansion]\label{lem:blocks}
Fix $K\ge1$. For $0\le j\le K$, let
$s_{j,N}=\ell_jN+v_{j,N}$ be integer endpoints with
$s_{0,N}=N+1$, $v_{0,N}=1$, and all other $v_{j,N}$ in a fixed bounded
set. For $1\le j\le K$, put
$B_{j,N}=\{s_{j,N},\ldots,s_{j-1,N}-1\}$. Uniformly over
these endpoints and $1\le j\le K$,
\begin{align}
 S_{j,N}:=\sum_{t\in B_{j,N}}o_t
 &=\lambda_j+\frac{\beta_{j,N}}N+O(N^{-2}),\label{eq:S}\\
 T_{j,N}:=\sum_{t\in B_{j,N}}o_t^2
 &=\frac{\delta_j}{N}+O(N^{-2}),\label{eq:T}
\end{align}
where
\begin{equation}\label{eq:beta}
 \beta_{j,N}
 =\frac{v_{j-1,N}-3/2}{\ell_{j-1}}
  -\frac{v_{j,N}-3/2}{\ell_j}.
\end{equation}
For every fixed integer $b\ge0$, with $\phi_b(x)=x^b/b!$,
\begin{equation}\label{eq:elementary}
 \esym_b(B_{j,N})=\phi_b(\lambda_j)
 +\frac1N\left(\beta_{j,N}\phi_b'(\lambda_j)
            -\frac{\delta_j}{2}\phi_b''(\lambda_j)\right)+O(N^{-2}).
\end{equation}
\end{lemma}
\begin{proof}
Write $\mathsf H_n=\sum_{j=1}^n j^{-1}$ and let $\gamma$ be Euler's constant.
For $s=\ell N+v$ with fixed $\ell>0$ and bounded $v$,
\[
 \mathsf H_{s-2}=\log N+\log\ell+\gamma
                 +\frac{v-3/2}{\ell N}+O(N^{-2}).
\]
Taking the difference at the block endpoints and using
$\log(\ell_{j-1}/\ell_j)=\lambda_j$ gives \eqref{eq:S}.
Comparison with the integral of $x^{-2}$ gives \eqref{eq:T} with
error $O(N^{-2})$. More generally,
$\sum_{t\in B_{j,N}}o_t^a=O(N^{1-a})$ for each fixed $a\ge3$.

For completeness, extract the coefficient of $z^b$ in
\[
 \prod_{t\in B_{j,N}}(1+zo_t)
 =\exp\left(zS_{j,N}-\frac{z^2}{2}T_{j,N}
       +\sum_{a\ge3}\frac{(-1)^{a+1}z^a}{a}
                   \sum_{t\in B_{j,N}}o_t^a\right).
\]
Only powers through degree $b$ matter. All terms involving a power sum
of degree at least three, or two factors $T_{j,N}$, are $O(N^{-2})$.
It follows that
\begin{equation}\label{eq:elementary-intermediate}
 \esym_b(B_{j,N})
 =\frac{S_{j,N}^b}{b!}
       -\frac{T_{j,N}S_{j,N}^{b-2}}{2(b-2)!}+O(N^{-2}),
\end{equation}
with the second term omitted for $b<2$. Substitution of
\eqref{eq:S}--\eqref{eq:T} proves \eqref{eq:elementary}. All estimates
are uniform because the displacements range over a bounded set and only
finitely many coefficients are extracted.
\end{proof}

\subsection{The derivative identity}

For $p\ge t\ge1$, define the truncated pattern polynomial
\[
 A_{p,t}(x_1,\ldots,x_t)
 =\sum_{\substack{\mathbf b\in\Xi_p\\b_p=\cdots=b_{t+1}=0}}
                      \prod_{j=1}^t\frac{x_j^{b_j}}{b_j!}.
\]
In particular, $A_{t,t}=H_t$ and $A_{k,k-1}=\Phi_k$. Sorting the
patterns by their $t$th coordinate gives
\begin{equation}\label{eq:A-recurrence}
 A_{p,t}=\frac{x_t^{p-t+1}}{(p-t+1)!}
       +\sum_{b=0}^{p-t}\frac{x_t^b}{b!}A_{p-b,t-1}\quad(t\ge2),
 \qquad A_{p,1}=\frac{x_1^p}{p!}.
\end{equation}
In the first term the active suffix in \eqref{eq:patterns} begins at $t$.
Otherwise $0\le b_t=b\le p-t$; keeping the lower coordinates and adjusting
the number of initial zero coordinates gives exactly a pattern counted
by $A_{p-b,t-1}$. This correspondence is reversible and preserves the
displayed weights.

\begin{lemma}[Derivative identity]\label{lem:derivatives}
If $p>t\ge1$ and $1\le j\le t$, then
\begin{equation}\label{eq:A-derivative}
 \partial_j A_{p,t}(\lambda_1,\ldots,\lambda_t)
       =A_{p-1,t}(\lambda_1,\ldots,\lambda_t).
\end{equation}
Consequently, for all $j<k$,
\begin{equation}\label{eq:unit-derivatives}
 \partial_j\Phi_k(\lambda_1,\ldots,\lambda_{k-1})=1.
\end{equation}
\end{lemma}
\begin{proof}
Induct on $t$. The case $t=1$ follows by differentiating $x_1^p/p!$.
For $j=t$, differentiating \eqref{eq:A-recurrence} gives the polynomial
identity $\partial_t A_{p,t}=A_{p-1,t}$. For $j<t$, the induction
hypothesis gives
\[
 \partial_j A_{p,t}(\lambda_1,\ldots,\lambda_t)
 =\sum_{b=0}^{p-t}\frac{\lambda_t^b}{b!}
          A_{p-b-1,t-1}(\lambda_1,\ldots,\lambda_{t-1}).
\]
The last term is $\lambda_t^{p-t}/(p-t)!$ because
$A_{t-1,t-1}(\lambda_1,\ldots,\lambda_{t-1})=1$.
The remaining terms are precisely the sum in the recurrence for
$A_{p-1,t}$. Finally, take $p=k$, $t=k-1$, and use
$A_{k-1,k-1}(\lambda_1,\ldots,\lambda_{k-1})=1$.
\end{proof}

Write
\begin{equation}\label{eq:h}
 h_{k,j}=\partial_{jj}\Phi_k(\lambda_1,\ldots,\lambda_{k-1})
 \quad(1\le j<k).
\end{equation}
These numbers are nonnegative rational numbers. Moreover,
\begin{equation}\label{eq:last-h}
 h_{k,k-1}=1\qquad(k\ge2).
\end{equation}
For $k=2$ this follows from $\Phi_2=x_1^2/2$; for $k\ge3$, use
$\Phi_k=x_{k-1}^2/2+x_{k-1}A_{k-1,k-2}+A_{k,k-2}$.

\section{Local threshold corrections}\label{sec:local}

For a function $F$ of $d$ variables, write
$F(\lam)=F(\lambda_1,\ldots,\lambda_d)$.

Set $V_{k,N}(r)=P_{k,N}(r+1)$ and
$\Delta_{k,N}(r)=V_{k,N}(r)-V_{k,N}(r-1)$. The following theorem gives
the correction at every fixed level.

\begin{theorem}[Local expansion]\label{thm:local}
For $k\ge1$, define
\begin{equation}\label{eq:c-definition}
 c_k=\frac12-\frac{\ell_k}{2}
          \left(1+\sum_{j<k}h_{k,j}\delta_j\right).
\end{equation}
Then $i_N^{(k)}-\ell_kN=O(1)$. For every fixed $R<\infty$,
\begin{equation}\label{eq:local}
 N^2\Delta_{k,N}(\ell_kN+u)
       =-\frac{u-c_k}{\ell_k}+O(N^{-1})
\end{equation}
uniformly over $|u|\le R$ for which $\ell_kN+u$ is an admissible integer
threshold. The constants in the error may depend on $k$ and $R$.
\end{theorem}
\begin{proof}
We prove bounded localization by induction. Assume that
$i_N^{(j)}=\ell_jN+v_{j,N}$ with bounded $v_{j,N}$ for $j<k$, and
consider $i=\ell_kN+v$ with bounded $v$. For $k=1$ there are no later
thresholds, so the same calculation starts the induction.
Applying Lemma~\ref{lem:blocks} to the factors defining
$\mathcal C_{k,N}$ yields
\begin{align}
 \mathcal C_{k,N}
 &=\Phi_k(\lambda_1,\ldots,\lambda_{k-1})\notag\\
 &\quad+\frac1N\sum_{j<k}
       \left(\beta_{j,N}\partial_j\Phi_k(\lam)
                    -\frac{\delta_j}{2}\partial_{jj}\Phi_k(\lam)\right)
       +O(N^{-2})\notag\\
 &=1-\lambda_k+\frac1N
       \left(\sum_{j<k}\beta_{j,N}
                    -\frac12\sum_{j<k}h_{k,j}\delta_j\right)
       +O(N^{-2}),\label{eq:C-expansion}
\end{align}
where \eqref{eq:lambda} and Lemma~\ref{lem:derivatives} give the last
line. The current block has odds sum
$\lambda_k+\beta_{k,N}/N+O(N^{-2})$.

The endpoint corrections telescope:
\begin{equation}\label{eq:telescoping}
 \sum_{j=1}^k\beta_{j,N}
  =\frac{1-3/2}{\ell_0}-\frac{v-3/2}{\ell_k}
  =-\frac12-\frac{v-3/2}{\ell_k}.
\end{equation}
In particular, every displacement of a later threshold has disappeared.
Substituting in \eqref{eq:difference} gives
\begin{equation}\label{eq:G-local}
 G_{k,N}(\ell_kN+v)
       =\frac{v-(c_k+1)}{N\ell_k}+O(N^{-2}).
\end{equation}

Choose a fixed $R>|c_k+1|+2$ and take integers within one unit of
$\ell_kN-R$ and $\ell_kN+R$. They are admissible for all sufficiently
large $N$, since $0<\ell_k<\ell_{k-1}$. Equation~\eqref{eq:G-local}
is negative at the first and positive at the second. Strict monotonicity
of $G_{k,N}$ places all maximizers between these integers, up to one
endpoint step. Hence $i_N^{(k)}-\ell_kN=O(1)$, closing the induction.

Finally, $r=\ell_kN+u$ corresponds to $i=r+1$, and
\[
 \Delta_{k,N}(r)=-\frac1N G_{k,N}(r+1).
\]
Putting $v=u+1$ in \eqref{eq:G-local} proves \eqref{eq:local}.
\end{proof}

\begin{remark}\label{rem:arbitrary-later}
The expansion \eqref{eq:local} itself does not require the later
thresholds to be optimal. The same proof applies if they are any
admissible integers $\ell_jN+O(1)$ in starting-index coordinates, with a
uniform error over bounded displacements. Optimality is used to identify
the maximizing threshold, not to obtain the cancellation.
\end{remark}

\subsection{Poisson counts and collisions}

\begin{proposition}[Factorial-moment interpretation]\label{prop:poisson}
Let $Z_1,\ldots,Z_k$ be independent Poisson variables of respective means
$\lambda_1,\ldots,\lambda_k$, and put
$\mathcal W_k=\{(Z_k,\ldots,Z_1)\in\Xi_k\}$.
Then $\Prob(\mathcal W_k)=\ell_k$, and for $j<k$,
\begin{equation}\label{eq:poisson-h}
 h_{k,j}=\frac1{\lambda_j^2}
             \E[Z_j(Z_j-1)\mid\mathcal W_k].
\end{equation}
Consequently,
\begin{equation}\label{eq:poisson-c}
 c_k=\frac12\Prob(\mathcal W_k^{\mathrm c})
       -\E\left[\mathbf1_{\mathcal W_k}
            \sum_{j<k}\frac{\delta_j}{\lambda_j^2}\binom{Z_j}{2}\right].
\end{equation}
\end{proposition}
\begin{proof}
Summing the Poisson masses over $\Xi_k$ gives
$\Prob(\mathcal W_k)=\e^{-L_k}H_k(\lam)=\ell_k$.
Multiplying those masses by $b_j(b_j-1)$ gives
\[
 \E[Z_j(Z_j-1)\mathbf1_{\mathcal W_k}]
      =\ell_k\lambda_j^2\partial_{jj}H_k(\lam)
      =\ell_k\lambda_j^2h_{k,j}.
\]
Division by $\ell_k$ proves \eqref{eq:poisson-h}; substitution into
\eqref{eq:c-definition} proves \eqref{eq:poisson-c}.
\end{proof}

Here is an exact finite-horizon explanation of the collision term.
For $2\le t\le N$, let $Y_t$ be independent Poisson variables of means
$o_t$. Since
\[
 \Prob(Y_t=1\mid Y_t\le1)=\frac{o_t}{1+o_t}=p_t,
\]
conditioning on $Y_t\le1$ for every $t$ gives precisely the joint law
of the record indicators. In a block $B$, conditional on
$\sum_{t\in B}Y_t=b$, the $b$ Poisson points are independently assigned
to candidate times with probabilities $o_t/S$, where
$S=\sum_{t\in B}o_t$. Thus, writing $\mathcal D$ for the event that
the assigned times are distinct,
\begin{equation}\label{eq:embedding}
 \esym_b(B)=\frac{S^b}{b!}\Prob(\mathcal D\mid\textstyle\sum_{t\in B}Y_t=b).
\end{equation}
These elementary Poisson conditioning properties are also discussed
in~\cite{LastPenrose2018}.

For the blocks of Lemma~\ref{lem:blocks}, inclusion--exclusion gives
\[
 \Prob(\mathcal D\mid\textstyle\sum_{t\in B}Y_t=b)
       =1-\binom b2\frac{T}{S^2}+O(N^{-2}),
 \qquad T=\sum_{t\in B}o_t^2.
\]
Indeed, a specified pair collides with probability $T/S^2$. Intersections
of two distinct pair-collision events are $O(N^{-2})$: disjoint pairs
give $(T/S^2)^2$, and overlapping pairs give
$\sum_{t\in B}(o_t/S)^3$. There are finitely many such pairs for fixed
$b$. Substituting in \eqref{eq:embedding} recovers
\eqref{eq:elementary-intermediate}. Only diagonal Hessian terms appear
because distinct blocks contain disjoint candidate times.

\section{The optimal value and bounded rounding}\label{sec:value}

The cancellation in \eqref{eq:telescoping} has a second consequence:
bounded rounding changes the value only at order $N^{-2}$.

\begin{theorem}[Uniform value expansion]\label{thm:value}
Fix $m\ge1$ and $R<\infty$. For any integer rejected-count vector
$\mathbf r=(r_m,\ldots,r_1)$ satisfying
$|r_k-\ell_kN|\le R$ for $1\le k\le m$, its threshold rule is
admissible for all sufficiently large $N$ and, with
$\mathbf1=(1,\ldots,1)\in\R^m$, has success probability
\begin{equation}\label{eq:value-uniform}
 P_N(\mathbf r+\mathbf1)
  =\sum_{k=1}^m\ell_k+\frac1N\sum_{k=1}^m c_k+O(N^{-2}),
\end{equation}
uniformly over all such vectors. In particular,
\begin{equation}\label{eq:value-optimal}
 P_N^*(m)=\sum_{k=1}^m\ell_k+\frac1N\sum_{k=1}^m c_k+O(N^{-2}).
\end{equation}
\end{theorem}
\begin{proof}
Since $0<\ell_m<\cdots<\ell_1<1$, bounded perturbations preserve
admissibility for large $N$. Write $r_k=\ell_kN+w_k$ and use starting
indices $i_k=r_k+1$, so the displacements in Lemma~\ref{lem:blocks}
are $v_k=w_k+1$.

All first derivatives of $H_k$ equal one at $\lam$: for $j<k$ this is
Lemma~\ref{lem:derivatives}, and $\partial_kH_k=1$ by
\eqref{eq:polynomials}. The diagonal second derivatives are $h_{k,j}$
for $j<k$ and zero for $j=k$. Therefore the same finite expansion as in
\eqref{eq:C-expansion} gives
\begin{equation}\label{eq:H-expansion}
 \mathcal H_{k,N}
 =1+\frac1N\left(\sum_{j=1}^k\beta_{j,N}
                 -\frac12\sum_{j<k}h_{k,j}\delta_j\right)+O(N^{-2}).
\end{equation}
By \eqref{eq:telescoping},
$\sum_{j=1}^k\beta_{j,N}=-1/2-(w_k-1/2)/\ell_k$.
Multiplying \eqref{eq:H-expansion} by $r_k/N=\ell_k+w_k/N$ consequently
yields
\begin{align*}
 \frac{r_k}{N}\mathcal H_{k,N}
 &=\ell_k+\frac1N
     \left(w_k-\frac{\ell_k}{2}-w_k+\frac12
                     -\frac{\ell_k}{2}\sum_{j<k}h_{k,j}\delta_j\right)
       +O(N^{-2})\\
 &=\ell_k+\frac{c_k}{N}+O(N^{-2}).
\end{align*}
Sum over $k$ in \eqref{eq:exact-value}. Uniformity follows from the
bounded displacements and fixed $m$. The optimal vector satisfies the
required boundedness by Theorem~\ref{thm:local}, proving
\eqref{eq:value-optimal}.
\end{proof}

\begin{corollary}[Loss from bounded rounding]\label{cor:regret}
Under the hypotheses of Theorem~\ref{thm:value}, there is a constant
$C_{m,R}>0$, depending only on $m$ and $R$, such that
\[
 0\le P_N^*(m)-P_N(\mathbf r+\mathbf1)\le C_{m,R}N^{-2}
\]
for all sufficiently large $N$. This applies, in particular, to
$r_k=\floor{\ell_kN}$ and to $r_k=\floor{\ell_kN+c_k}$.
\end{corollary}
\begin{proof}
Nonnegativity follows from optimality. Subtract
\eqref{eq:value-uniform} from \eqref{eq:value-optimal}; their constant
and $N^{-1}$ terms are identical. Both displayed rounding rules have
bounded displacements for fixed $m$.
\end{proof}

The corollary separates two levels of approximation. A boundedly rounded
vector has the optimal first finite-horizon value correction, while the
sign of the local increment is needed to identify the exact maximizing
integer. We now address the latter question.

\section{Integer rounding and arithmetic consequences}\label{sec:rounding}

\subsection{The rounding law}

Write $\fract{x}=x-\floor{x}$ for the fractional part.
The correction need not be positive. Set
\begin{equation}\label{eq:ad}
 a_k=\floor{c_k},\qquad d_k=c_k-a_k.
\end{equation}
We first exclude an integer value of $c_k$.

\begin{proposition}[Nonintegral corrections]\label{prop:transcendence}
For every $k\ge1$, both $\ell_k$ and $c_k$ are transcendental.
In particular, $0<d_k<1$.
\end{proposition}
\begin{proof}
The numbers $L_j$ and $h_{k,j}$ are rational. Since $L_k>0$,
transcendence of $\ell_k=\e^{-L_k}$ follows from the
Lindemann--Weierstrass theorem~\cite{Baker1975}. Rewrite
\eqref{eq:c-definition} as
\begin{equation}\label{eq:c-exponentials}
 2c_k=1-\e^{-L_k}
         \left(1+\sum_{j<k}h_{k,j}(\e^{L_j}-\e^{L_{j-1}})\right).
\end{equation}
For $k=1$ this is $2c_1=1-\e^{-1}$. For $k\ge2$, after collecting
equal exponents the term $\e^{-(L_k-L_{k-1})}=\e^{-\lambda_k}$ has
coefficient $-1$, by \eqref{eq:last-h}. Its exponent is distinct from
all the others because $L_0<\cdots<L_k$. If $c_k$ were algebraic,
\eqref{eq:c-exponentials} would give a nontrivial algebraic linear
relation among exponentials of distinct algebraic numbers, contrary
to Lindemann--Weierstrass.
\end{proof}

\begin{theorem}[Asymptotic rounding]\label{thm:rounding}
For each fixed $k$, all sufficiently large $N$ satisfy
\begin{equation}\label{eq:two-point}
 r_N^{(k)}\in
 \left\{\floor{\ell_kN}+a_k,\ \floor{\ell_kN}+a_k+1\right\}.
\end{equation}
There is a constant $C_k>0$ such that, whenever
\begin{equation}\label{eq:transition}
 \left|\fract{\ell_kN}-(1-d_k)\right|>\frac{C_k}{N},
\end{equation}
the $k$th optimum with the later optimal thresholds fixed is unique and
\begin{equation}\label{eq:rounding}
 r_N^{(k)}=\floor{\ell_kN+c_k}.
\end{equation}
\end{theorem}
\begin{proof}
Let $m_N=\floor{\ell_kN}$, $\theta_N=\fract{\ell_kN}$, and
$r_0=m_N+a_k$. The offsets of $r_0$ and $r_0+2$ are
$a_k-\theta_N$ and $a_k+2-\theta_N$. The first is uniformly below
$c_k$ by at least $d_k$, and the second is uniformly above $c_k$
by at least $1-d_k$. Theorem~\ref{thm:local} therefore gives
$\Delta_{k,N}(r_0)>0$ and $\Delta_{k,N}(r_0+2)<0$ for large $N$.
Unimodality proves \eqref{eq:two-point}.

For the remaining comparison, the offset of $r_0+1$ is
$a_k+1-\theta_N$, and \eqref{eq:local} gives
\[
 N^2\Delta_{k,N}(r_0+1)
      =\frac{\theta_N-(1-d_k)}{\ell_k}+O(N^{-1}).
\]
Under \eqref{eq:transition}, its sign is negative
when $\theta_N<1-d_k$ and positive when $\theta_N>1-d_k$.
These two cases select $r_0$ and $r_0+1$, respectively, and give
\eqref{eq:rounding} with a strict maximum.
\end{proof}

The transition band is an error-controlled limitation of the formula.
The theorem does not assert that the exceptional set is finite.
Resolving individual horizons inside that band requires a higher-order
expansion or the exact difference \eqref{eq:difference}.

\subsection{Approximating integers and their frequencies}

\begin{corollary}[Accurate integer approximations]\label{cor:approximations}
Fix $k$. Let $q_j$ be positive integers tending to infinity and let
$p_j$ be integers. If $p_j-\ell_kq_j\to0$, then
\begin{equation}\label{eq:arithmetic}
 r_{q_j}^{(k)}=p_j+a_k
\end{equation}
for all sufficiently large $j$. In particular this holds when $p_j/q_j$
are the continued-fraction convergents of $\ell_k$.
\end{corollary}
\begin{proof}
At horizon $q_j$, the offsets of $p_j+a_k$ and $p_j+a_k+1$ converge
to $a_k$ and $a_k+1$. Since $a_k<c_k<a_k+1$, the local increments
are eventually positive at the first integer and negative at the next.
Unimodality identifies the unique optimum. Convergents satisfy
$|p_j-\ell_kq_j|<1/q_j$, so they meet the hypothesis.
\end{proof}

The natural density of $A\subseteq\N$ is
$\lim_{M\to\infty}\#(A\cap\{1,\ldots,M\})/M$, when this limit exists.

\begin{corollary}[Joint rounding frequencies]\label{cor:frequencies}
Fix $m\ge1$. For every $\boldsymbol\varepsilon=(\varepsilon_1,\ldots,
\varepsilon_m)\in\{0,1\}^m$, the set of horizons satisfying
\[
 r_N^{(k)}=\floor{\ell_kN}+a_k+\varepsilon_k
 \qquad(1\le k\le m)
\]
has natural density
\begin{equation}\label{eq:joint-frequency}
 \prod_{k=1}^m d_k^{\varepsilon_k}(1-d_k)^{1-\varepsilon_k}.
\end{equation}
In particular, the two choices for coordinate $k$ have respective
frequencies $1-d_k$ and $d_k$.
\end{corollary}
\begin{proof}
The numbers $0,-L_1,\ldots,-L_m$ are distinct algebraic numbers.
Lindemann--Weierstrass therefore implies that
$1,\ell_1,\ldots,\ell_m$ are linearly independent over the rationals.
The Kronecker--Weyl equidistribution theorem
\cite{KuipersNiederreiter1974} now gives equidistribution of
\[
 \bigl(\fract{\ell_1N},\ldots,\fract{\ell_mN}\bigr)
 \quad\text{on }[0,1)^m.
\]
For each $k$, the exceptional band in \eqref{eq:transition} has density
zero: for each $\varepsilon>0$, it is eventually contained in the fixed
interval of radius $\varepsilon$ about $1-d_k$, whose visit frequency
is at most $2\varepsilon$. Letting $\varepsilon\downarrow0$ proves
this assertion, and the union of the $m$ bands also has density zero.
Outside that union, Theorem~\ref{thm:rounding} associates each binary
vector with a rectangular box whose $k$th side has length $1-d_k$
if $\varepsilon_k=0$ and $d_k$ if $\varepsilon_k=1$. Its volume is
\eqref{eq:joint-frequency}.
\end{proof}

The factorization concerns frequencies as the integer horizon varies.
It does not assert independence of the finite-horizon optimization
problems. It records the combined effect of the local cancellation
and rational independence of the limiting proportions.

\section{Examples and exact finite certification}\label{sec:examples}

\subsection{The first three corrections}

For one selection, $\Phi_1=0$, $\lambda_1=1$, and
\begin{equation}\label{eq:c1}
 \ell_1=\e^{-1},\qquad c_1=\frac{1-\e^{-1}}2.
\end{equation}
For two selections, $\Phi_2=x_1^2/2$, so $\lambda_2=1/2$ and
$h_{2,1}=1$. Thus
\begin{equation}\label{eq:c2}
 \ell_2=\e^{-3/2},\qquad c_2=\frac{1-\e^{-1/2}}2.
\end{equation}
For three selections,
\[
 \Phi_3(x_1,x_2)=\frac{x_2^2}{2}
                  +\frac{x_2x_1^2}{2}+\frac{x_1^3}{6},
 \qquad \lambda_3=\frac{11}{24},
 \qquad (h_{3,1},h_{3,2})=\left(\frac32,1\right).
\]
Consequently
\begin{equation}\label{eq:c3}
 \ell_3=\e^{-47/24},\qquad
 c_3=\frac12-\frac14
        \left(2\e^{-11/24}+\e^{-23/24}-\e^{-47/24}\right).
\end{equation}
For instance, the three-selection optimal value is
\[
 P_N^*(3)=\e^{-1}+\e^{-3/2}+\e^{-47/24}
                +\frac{c_1+c_2+c_3}{N}+O(N^{-2}).
\]
The individual integer thresholds oscillate, but their displacements do
not enter this coefficient.

\subsection{The first negative correction}

The recurrence \eqref{eq:A-recurrence} computes $\lambda_k$ and
$h_{k,j}$ with exact rational arithmetic, without enumerating $\Xi_k$.
Table~\ref{tab:constants} records the first twelve levels.

\begin{table}[t]
\centering
\small
\caption{Leading parameters and threshold corrections. Entries are rounded
to twelve decimal places; the sign certification uses exact rational
intervals, not the displayed decimals.}\label{tab:constants}
\begin{tabular}{@{}rccc@{}}
\toprule
$k$&$\lambda_k$&$\ell_k$&$c_k$\\
\midrule
1&1.000000000000&0.367879441171& 0.316060279414\\
2&0.500000000000&0.223130160148& 0.196734670144\\
3&0.458333333333&0.141093380701& 0.123222120863\\
4&0.438368055556&0.091017690625& 0.077139226712\\
5&0.426267496745&0.059429241914& 0.047983259454\\
6&0.418025602017&0.039124966449& 0.029456969884\\
7&0.411997951369&0.025913468166& 0.017676677613\\
8&0.407371710247&0.017242739041& 0.010204982483\\
9&0.403694517656&0.011515530609& 0.005494627886\\
10&0.400692808363&0.007713745010& 0.002555937036\\
11&0.398190596376&0.005180042223& 0.000752738780\\
12&0.396069113648&0.003485962165&$-0.000325106395$\\
\bottomrule
\end{tabular}
\end{table}

\begin{proposition}[Certified finite signs]\label{prop:signs}
The constants satisfy $0<c_k<1/2$ for $1\le k\le11$ and
$-1<c_{12}<0$. More precisely,
\begin{equation}\label{eq:sign-margins}
 \frac{752738779}{10^{12}}<c_{11}<\frac{752738781}{10^{12}},
 \qquad
 -\frac{325106396}{10^{12}}<c_{12}< -\frac{325106394}{10^{12}}.
\end{equation}
\end{proposition}
\begin{proof}[Exact-arithmetic certificate]
Differentiate \eqref{eq:A-recurrence} twice in each variable and
evaluate recursively at the rational parameters obtained from
\eqref{eq:lambda}. This determines every required $h_{k,j}$ exactly.
For a nonnegative rational $x$, set
\[
 E_M(x)=\sum_{a=0}^M\frac{x^a}{a!},\qquad
 R_M(x)=\frac{x^{M+1}}{(M+1)!}
                         \frac1{1-x/(M+2)}.
\]
If $x<M+2$, comparison of successive positive Taylor terms gives
\begin{equation}\label{eq:exp-certificate}
 E_M(x)\le\e^x\le E_M(x)+R_M(x).
\end{equation}
Use $M=100$ and $x=L_j$, $0\le j\le12$. All these inputs are
rational and less than $102$. Let $[E_j^-,E_j^+]$ be the resulting
interval for $\e^{L_j}$, with $E_0^-=E_0^+=1$. Then
\[
 D_j^-=E_j^--E_{j-1}^+,
 \qquad D_j^+=E_j^+-E_{j-1}^-
\]
enclose $\e^{L_j}-\e^{L_{j-1}}$; each $D_j^->0$ is verified by a
rational comparison. Since $h_{k,j}\ge0$, the quantities
\[
 B_k^-=1+\sum_{j<k}h_{k,j}D_j^-,\qquad
 B_k^+=1+\sum_{j<k}h_{k,j}D_j^+
\]
are positive bounds for the bracket in \eqref{eq:c-exponentials}.
Therefore
\begin{equation}\label{eq:c-certificate}
 \frac12\left(1-\frac{B_k^+}{E_k^-}\right)
 \le c_k\le
 \frac12\left(1-\frac{B_k^-}{E_k^+}\right).
\end{equation}
Evaluation of these rational endpoints gives the stated signs and
\eqref{eq:sign-margins}. The accompanying
\texttt{compute\_constants.py} implements exactly these finite recurrences
using integer fractions and checks each claimed inequality. It also checks
$A_{k,k}(\lam)=1$, \eqref{eq:unit-derivatives}, and
\eqref{eq:last-h} at every computed level. No floating-point comparison
is used in the certificate.
\end{proof}

Only Proposition~\ref{prop:signs} uses this finite certificate; the
threshold, value, and rounding theorems hold analytically for every
fixed $k$ and do not depend on computed signs.

\begin{corollary}[The first shifted convergent threshold]\label{cor:first-shift}
For $1\le k\le11$, sufficiently large convergents $p_j/q_j$ of
$\ell_k$ satisfy $r_{q_j}^{(k)}=p_j$. For convergents of $\ell_{12}$,
\[
 r_{q_j}^{(12)}=p_j-1
\]
for all sufficiently large $j$.
\end{corollary}
\begin{proof}
Proposition~\ref{prop:signs} gives $a_k=0$ for $k\le11$ and
$a_{12}=-1$. Apply Corollary~\ref{cor:approximations}.
\end{proof}

At level twelve the two eventual rejected-count thresholds are
$\floor{\ell_{12}N}-1$ and $\floor{\ell_{12}N}$.
Their respective densities are $-c_{12}$ and $1+c_{12}$, approximately
$0.000325106395$ and $0.999674893605$. Thus the negative correction
changes the convergent identity even though the lower of the two
possible thresholds occurs at a small proportion of horizons.
These are limiting densities; no rate of convergence for their
finite-horizon frequencies is claimed. No sign assertion for $k\ge13$
is needed here.

\subsection{Certified thresholds at convergent denominators}

The small magnitude of $c_{12}$ requires care when interpreting
finite-horizon checks. Table~\ref{tab:convergent-thresholds} gives certified
exact thresholds at six consecutive continued-fraction convergents of
$\ell_{12}$, together with a comparison of the finite-horizon gain and its
limiting correction. The continued-fraction status is certified by
extracting common partial quotients from the two rational endpoints for
$\ell_{12}$ obtained using \eqref{eq:exp-certificate}.

For completeness, let $U_k(t)$ be $N$ times the optimal success probability
when only record times $t,\ldots,N$ remain and $k$ selections are available.
Set $U_0(t)=0$ and $U_k(N+1)=0$. Conditional on a record at time $t$,
selecting it wins either when there is no later record, with probability
$t/N$, or when the remaining $k-1$ selections win later. These events are
disjoint. Independence of the record indicators therefore gives
\begin{equation}\label{eq:certified-bellman}
 U_k(t)=\max\left\{U_k(t+1),\,
 \frac{(t-1)U_k(t+1)+t+U_{k-1}(t+1)}{t}\right\}.
\end{equation}
The record is accepted precisely when
$g_k(t):=t+U_{k-1}(t+1)-U_k(t+1)>0$, with rejection at a tie.
At horizon $N=q$, define the finite-horizon residual $\rho_q$ by
\begin{equation}\label{eq:finite-gain-residual}
 g_{12}(p)=(p-\ell_{12}q)-c_{12}+\rho_q.
\end{equation}
Here $g_{12}(p)$ is evaluated by the finite Bellman recursion; $\rho_q$
is measured from this exact gain, rather than estimated by an asymptotic
remainder bound.

\suppressfloats[t]
\begin{table}[t]
\centering
\caption{Exact level-12 rejected-count thresholds at six consecutive
convergents $p/q$ of $\ell_{12}$. The last column displays the certified
gain residual from \eqref{eq:finite-gain-residual}, divided by $|c_{12}|$
and rounded to three decimals. The residual dominates in the first four
rows; the correction dominates in the last two, as quantified below.}
\label{tab:convergent-thresholds}
\begin{tabular}{@{}rrrrr@{}}
\toprule
$p$ & $q$ & $r_q^{(12)}$ & $r_q^{(12)}-p$ & $\rho_q/|c_{12}|$\\
\midrule
7 & $2\,008$ & 6 & $-1$ & 42.658\\
15 & $4\,303$ & 14 & $-1$ & 16.277\\
22 & $6\,311$ & 21 & $-1$ & 11.580\\
37 & $10\,614$ & 36 & $-1$ & 2.507\\
\midrule
$1\,428$ & $409\,643$ & $1\,427$ & $-1$ & 0.071\\
$4\,321$ & $1\,239\,543$ & $4\,320$ & $-1$ & 0.037\\
\bottomrule
\end{tabular}
\end{table}

The first four rows satisfy
$\rho_q>|p-\ell_{12}q|+|c_{12}|$: their finite-horizon residual alone
overwhelms the correction. In the last two rows,
\[
 |p-\ell_{12}q|+|\rho_q|<|c_{12}|,
\]
so the negative correction determines the sign of the gain even after
both the approximation error and the finite-horizon residual are included.
These inequalities are certified using rational enclosures for all terms
in \eqref{eq:finite-gain-residual}.

A contrasting calculation at level eleven shows why the distinction
matters. Although $c_{11}>0$, the convergent denominators
$q=3861,4054,7915$ give $r_q^{(11)}-p=-1$, whereas
$q=27799,146910,321619$ give $r_q^{(11)}-p=0$.
Thus an observed shift of $-1$ at a small denominator can be a transient
effect even when the eventual convergent shift is zero. The level-eleven
cases are certified by the same finite Bellman recursion. The exact sign
certificate in Proposition~\ref{prop:signs}, together with
Corollary~\ref{cor:first-shift}, establishes the eventual distinction
between levels eleven and twelve.

The program \nolinkurl{verify_convergent_thresholds.py} evaluates
\eqref{eq:certified-bellman} for every $1\le k\le m$ and $1\le t\le N$,
with $m=12$ for Table~\ref{tab:convergent-thresholds} and $m=11$ for the
contrast above. It also records two analogous level-ten cases.
It propagates rational intervals with denominator $2^{128}$ and rounds
outwards after each division. The update is nondecreasing in both future
values, so endpoint evaluation preserves rigorous enclosures.
Every gain has a certified strict sign. In particular,
$g_m(r)<0<g_m(r+1)$ at each reported $r$, and the full sequence
of decisions verifies its threshold form.
For $N=2008,4303,6311$ and $m=12$, a second implementation using unrounded
rational numbers gives the same twelve thresholds and values within the
certified intervals. The supplementary JSON file records the exact
integer interval endpoints, including those for the gain decomposition
and the contrasting levels.

\section{Discussion}\label{sec:discussion}

The first finite-horizon correction is governed by two features of the
winning-pattern polynomials: their first derivatives all equal one at
the limiting parameters, while their diagonal Hessian entries measure
within-block collisions. The first property removes the later threshold
displacements from the local expansion. The second supplies the explicit
constant $c_k$. The same calculation explains why all boundedly rounded
threshold vectors share the coefficient $\sum_{k=1}^m c_k$ in the success
probability.

Three extensions remain natural. A further local expansion would resolve
more horizons inside the shrinking transition band; at that order the
later threshold displacements may reappear. Extending the analysis to
other Bernoulli arrays requires appropriate uniform block expansions and
control of the corresponding payoff prefactor. Finally, allowing the
number of selection opportunities to grow with $N$ requires uniform
estimates in the level $k$.

\section*{Statements and declarations}
\noindent\textbf{Funding.} No funds, grants, or other support were received
during the preparation of this manuscript.

\smallskip
\noindent\textbf{Competing interests.} The authors declare no competing
interests.

\smallskip
\begingroup\emergencystretch=2em
\noindent\textbf{Data and code availability.}
The supplementary source package contains the exact-arithmetic program
\texttt{compute\_constants.py}, the independent verification program
\nolinkurl{verify_finite_horizon.py}, and the finite-horizon certificate
\nolinkurl{verify_convergent_thresholds.py}, together with their outputs.
The finite-horizon checks supplement the analytic asymptotic proofs.
\par\endgroup

\smallskip
\noindent\textbf{Use of generative AI.} OpenAI's ChatGPT was used to assist with computational verification, and manuscript preparation.

\bigskip
\begingroup\small
\noindent\textsc{Jos\'e Mar\'ia Grau Ribas}\\
Departamento de Matem\'aticas, Universidad de Oviedo, Avda. Calvo Sotelo s/n,\\
33007 Oviedo, Spain\\
\textit{E-mail address:} \href{mailto:grau@uniovi.es}{\texttt{grau@uniovi.es}}

\medskip
\noindent\textsc{Ra\'ul S\'anchez Gal\'an}\\
4i Intelligent Insights, Tecnoincubadora Marie Curie, PCT Cartuja,\\
41092 Sevilla, Spain\\
\textit{E-mail address:} \href{mailto:r.sanchez@4i.ai}{\texttt{r.sanchez@4i.ai}}
\par\endgroup

\end{document}